\magnification 1200 \hsize=140mm \vsize=200mm \hoffset=-4mm \voffset=-1mm
\pretolerance=500 \tolerance=1000 \brokenpenalty=5000
\frenchspacing
\def\tvi{\vrule height 12pt depth 5pt width 0pt}
\font\tenbb=msbm10
\font\sevenbb=msbm7
\font\fivebb=msbm5
\def\tv{\tvi\vrule}
\def\hfq{\hfill\quad}
\newfam\bbfam\textfont\bbfam=\tenbb
\scriptfont\bbfam=\sevenbb
\scriptscriptfont\bbfam=\fivebb

\def\abs#1{\left\vert#1\right\vert}

\def\cqfd{\unskip\kern 6pt\penalty 500
\raise -2pt\hbox{\vrule\vbox to 10pt{\hrule width 4pt\vfill
\hrule}\vrule}\par}

\def\frac#1#2{{#1\over#2}}
\def\tableau#1{\vcenter{\offinterlineskip\halign{%
\tv\hfq$\vcenter{\vskip2pt\hbox{$\displaystyle{##}$}\vskip 2pt}$\hfq\tv&&\hfq$\vcenter{\vskip2pt\hbox{$\displaystyle{##}$}\vskip2pt}$\hfq\tv\cr\noalign{\hrule}#1}}}
\def\crh{\cr\noalign{\hrule}}
\def\bb{\fam\bbfam}
\def\date{{\the\day}\ \ifcase\month\or Janvier\or F\'evrier\or Mars\or Avril
\or Mai\or Juin\or Juillet\or Ao\^ut\or Septembre \or Octobre\or Novembre\or D\'ecembre\fi\ {\the\year}}

  \hsize=140mm \vsize=200mm \hoffset=-4mm \voffset=-1mm
\pretolerance=500 \tolerance=1000 \brokenpenalty=5000
\frenchspacing
\def\tvi{\vrule height 12pt depth 5pt width 0pt}
\font\tenbb=msbm10
\font\sevenbb=msbm7
\font\fivebb=msbm5
\def\tv{\tvi\vrule}
\def\hfq{\hfill\quad}
\newfam\bbfam\textfont\bbfam=\tenbb
\scriptfont\bbfam=\sevenbb
\scriptscriptfont\bbfam=\fivebb

\def\abs#1{\left\vert#1\right\vert}

\def\cqfd{\unskip\kern 6pt\penalty 500
\raise -2pt\hbox{\vrule\vbox to 10pt{\hrule width 4pt\vfill
\hrule}\vrule}\par}

\def\frac#1#2{{#1\over#2}}
\def\tableau#1{\vcenter{\offinterlineskip\halign{%
\tv\hfq$\vcenter{\vskip2pt\hbox{$\displaystyle{##}$}\vskip 2pt}$\hfq\tv&&\hfq$\vcenter{\vskip2pt\hbox{$\displaystyle{##}$}\vskip2pt}$\hfq\tv\cr\noalign{\hrule}#1}}}
\def\crh{\cr\noalign{\hrule}}
\def\bb{\fam\bbfam}

\def\virg{\raise 2pt\hbox{,}}

\let\ds=\displaystyle
\date

\def\abs#1{\left\vert #1\right\vert}

\def\cqfd{\kern 6pt\hbox{\vrule\vbox to 6pt{\hrule width
                        6pt\vfil\hrule}\vrule}}

\def\a{\alpha}

\def\tend{\longrightarrow}
\def\z{\zeta}

\centerline{\bf Elementary methods for evaluating Jordan's sums }
$$\sum_{n\geq 1}\Big(1+\frac13+\cdots+\frac1{2n-1}\Big){1\over n^{2a}}\,\,{\rm and}\,\,
\sum_{n\geq 1}\Big(1+\frac13+\cdots+\frac1{2n-1}\Big){1\over (2n-1)^{2a}} $$
\centerline{\bf and analogous Euler's type sums and for setting a $\sigma$-sum theorem } \vskip 0,3cm
 \vskip 1cm
\centerline
{\bf G. Bastien}\vskip 0,6cm
\centerline {Institut Math\'ematique de Jussieu et CNRS, }
\centerline{ Universit\'e Pierre et Marie Curie, 4 place Jussieu 75252 PARIS CEDEX 5}

\centerline{email adress: bastien@math.jussieu.fr}

\vskip 0,3cm
\noindent{\bf Abstract}
\vskip 0,3cm
Our aim is first  to calculate the following sums : $$\sum_{n\geq 1}\Big(1+\frac13+\cdots+\frac1{2n-1}\Big){1\over n^{2a}}\,\,{\rm and}\,\,
\sum_{n\geq 1}\Big(1+\frac13+\cdots+\frac1{2n-1}\Big){1\over (2n-1)^{2a}}\eqno (J) $$ when  $a$ is an   integer $\geq 1$, by mean of double sum
methods or   integral representation.  The values of sum in $(J)$ are due  P.F. Jordan, in  {\sl Infinite sums
of psi  functions}, Bulletin of American Mathematical Society  (79) 4, 1973 . We think that we give here the first simple  an elementary proof of $(J)$ and of
other formulas which  are deduced in the litterature from Jordan ones. As a consequence of our calculations, we find an expression for  the sums
$\ds \sum_{n\geq 1}\Big(1+\frac12+\cdots+\frac1n\Big){1\over(2n+1)^{2a}}\cdot$
Moreover, in the last section, we find some relations between the sums $\ds \sigma(s,t):= \sum_{n\geq 1}\Big(1+\frac1{3^t}+\cdots+\frac1{{(2n-1)}^t}\Big){1\over
n^{s}}$ (with
$s$ and
$t$ integers, $s\geq 2$), which give the sums $\ds \sum_{n\geq 1}\Big(1+\frac12+\cdots+\frac1n\Big){1\over(2n+1)^{2a+1}}\cdot$ Finally, we prove a new
additive relation on the $\sigma(s,t)$ of same weight $s+t$ which is a $``\sigma$-sum theorem" analogous to this one involving the classical Euler's sums
$\z(s,t)$. 


\vskip 0,3cm
\noindent{\bf Introduction and notations}\vskip 0,3cm

We define harmonic numbers $H_n$ and semi-harmonic numbers $S_n$ by $$H_0=0, H_n=H_{n-1}+\frac1n \,\,(n\geq 1),\,\,\,\,\,S_0=0,S_n=S_{n-1}+\frac 1{2n-1}
\,\,(n\geq 1)\eqno(1.1)$$ and we put
$$J(b) :=\sum_{n\geq 1}{S_n\over n^{b}}\virg\,\,\,\bar J(b): =\sum_{n\geq 1}{S_n\over(2 n-1)^{b}}\cdot\eqno(1.2)$$These series converge when  $b>1.$
and we denote them by {\bf ``Jordan's sums"}

We will give expression of  $J(2a)$ and $\bar J(2a)$ when  $a$ is an integer $\geq 1$. In the present paper, we propose an elementary and autonomous  proof of
their value. In [5], R.Sitaramachandrarao uses the values of $J(2a) $ and $\bar J(2a)$ given by Jordan in [4] to evaluate some series analogous to the so called
Euler's sums  $\ds
\sum_{n\geq 1}\Big(1+\frac12+\cdots+\frac1{n}\Big){1\over n^{b}}\virg$ as, for example  $\displaystyle\sum_{n=1}^\infty (-1)^{n-1}{H_n\over n^{2a}}\cdot$ In this paper, we
inverse the order of the proof, deducing for instance $J(2a)$ from the previous sum that we prove directly.

We  use   $\ds \lambda (s)=\sum _{n=1}^\infty {1\over (2n-1)^s}=\bigg (1-{1\over 2^s}\bigg)\z(s)$, where $\displaystyle \z(s)=\sum _{n=1}^\infty {1\over n^s}$ for $s>1$. 
We use also    the classical Euler star  sum$$\zeta^*(b,1)=\sum _{n\geq 1}{H_n\over n^b}=\Big(1+\frac
b2\Big)\z(b+1)-\frac12\sum_{j=2}^{b-1}\z(j)\z(b+1-j)\eqno (1.3)$$ when the integer $b$ is greater than 1. This  expression has a simpler form for even $b$ ( say $2a$)
$$\zeta^*(2a,1)=(1+a)\z(2a+1)-\sum_{j=1}^{a-1}\z(2j)\z(2a+1-2j).\eqno (1.4)$$
We denote by $\widetilde\zeta^*(b,1)$ the alternating  corresponding sum, that is:$$\widetilde\zeta^*(b,1)=\sum _{n\geq 1}(-1)^{n-1}{H_n\over n^b}\cdot$$

 Besides the formulas giving $J(2a)$ and $\bar J(2a)$, we obtain close formulas for the sums   $\ds \sum_{n\geq 1}{H_n\over (2n+1)^q}$  when $q$ in
an integer $\geq 2$.  The result when $q$ is odd is obtained by searching linear relations between  the sums $\ds \sigma(s,t):=
\sum_{n\geq 1}\Big(1+\frac1{3^t}+\cdots+\frac1{{(2n-1)}^t}\Big){1\over n^{s}}$ (with
$s$ and
$t$ integers, $s\geq 2$). We give also  explicit values of some of these  $\sigma (s,t)$, which we call $\sigma $-Euler sums in term
of Riemann's series or analogous ones and establish a sum theorem for all $\sigma$'  series  of same weight.   Finally, we sketch out a calculation of them
when
$s+t$ is odd.\vskip 0,8 cm

\noindent{\bf I. The results :

Principal formulas for Jordan's and analogous sums and $\sigma$-Euler sums}\vskip 0,5cm

\noindent{\bf A.\hskip 4 mm Jordan's sums}

$$\tableau{J(2a)=\sum_{n\geq 1}{S_n\over
n^{2a}}={2^{2a+1}-1\over 4}\z(2a+1)-\frac12\sum_{j=1}^{a-1}(2^{2j+1}-1)\z(2j+1)\z(2a-2j)\crh}\,\eqno(a)$$
This formula becomes, by using $\lambda$'s series: 
$$\tableau{J(2a)=\sum_{n\geq 1}{S_n\over
n^{2a}}=2^{2a-1}\lambda(2a+1)-\sum_{j=1}^{a-1}2^{2j}\lambda(2j+1)\z(2a-2j).\crh}\eqno(a')$$

 $$\tableau{\bar J(2a)=\sum_{n\geq 1}{S_n\over
(2n-1)^{2a}}=\lambda(2a)\ln2+\frac12\lambda(2a+1)-\sum_{j=1}^{a-1}{1\over 2^{2j+1}}\lambda(2a-2j)\z(2j+1).\crh}\eqno(b)$$

These are, with  slightly different notations, the formulas given by  Jordan in [4].\vskip 2mm
\noindent{\bf B.\hskip 2 mm Sums involving $H_{n}\,'s$}

We get two classes of closed formulas :$$\tableau{\sum_{n\geq 1}{H_{n}\over (2n+1
)^{2a}}=-2\lambda(2a)\ln2+2a\lambda(2a+1)-2\sum_{j=1}^{a-1}\lambda(2j)\lambda(2a+1-2j)\crh}\,.\eqno(c)$$
$$\tableau{\sum_{n\geq 1}{H_n\over(2n+1)^{2a-1}}=-2\lambda(2a-1)\ln
2+\Big(a-\frac12\Big)\lambda(2a)-\sum_{q=1}^{a-2}\lambda(2q+1)\lambda(2a-2q-1)\crh}\,\eqno(d)$$ 
When $a=2b$ it becomes: 
$$\tableau{\sum_{n\geq 1}{H_n\over(2n+1)^{4b-1}}=-2\lambda(4b-1)\ln
2+\Big(2b-\frac12\Big)\lambda(4b)-2\sum_{q=1}^{b-1}\lambda(2q+1)\lambda(4b-2q-1)\crh},\eqno(e)$$ and when $a=2b+1$ :
$$\tableau{\eqalign{\sum_{n\geq 1}{H_n\over(2n+1)^{4b+1}}=&-2\lambda(4b+1)\ln
2+\Big(2b+\frac12\Big)\lambda(4b+2)\cr&-\lambda^2(2b+1)-2\sum_{q=1}^{b-1}\lambda(2q+1)\lambda(4b-2q+1)}\crh}.\eqno(e')$$

\noindent {\bf C.\hskip 2 mm Some close formulas for $\sigma$-Euler sums and sum theorem}

Clearly, by definition, $\sigma(2a,1)=J(2a)$. We prove  the two  relations :  $$\tableau{\sigma(2,2a-1)=2a(2a-1)\lambda(2a+1)-8\sum_{j=1}^{a-1}j\lambda(2a-2j)\lambda(2j+1)\crh},\eqno(f)$$

$$\tableau{\eqalign{\sigma(2a-1,2)=-a2^{2a-1}\lambda(2a+1)+{2^{2a-1}(2a+1)\over
3}\lambda(2)\lambda(2a-1)\cr+\sum_{j=1}^{a-2}j2^{2j}\lambda(2j+1)\z(2a-2j)}\crh}.\eqno(g)$$
 
\proclaim Theorem. [of the $\sigma $-sum] The sum of all the $\sigma$'s of same weight is calculable: for $w\geq 3$, we have:
$$\tableau{\ds\sum_{i=1}^{w-2}\sigma(w-i,i)=(w-1)\lambda(w)\crh}.\eqno(h)$$

\vskip 0,3cm 

\noindent{\bf II. The  first result: evaluation of sums involving the $H_{2n}$}'s\vskip 0,5cm

The principal objective of this section is to evaluate the new sum $Z(2a):=\displaystyle \sum_{n\geq 1}{H_{2n}\over
n^{2a}}\cdot$ We begin by the  sum  $\displaystyle\widetilde\z^*(2a,1).$

\noindent{\bf II.1 Calculation of $\displaystyle\widetilde\z^*(2a,1):=\sum_{n\geq 1}(-1)^{n-1}{H_{n}\over
n^{2a}}$}

By rational decomposition, we get  $\ds {H_n\over n}=\sum_{q\geq 1}{1\over q(q+n)}\cdot$ So for $a\geq 1$,
by absolute convergence$$\widetilde\z^*(2a,1)=\sum_{n\geq 1 }{(-1)^{n-1}\over n^{2a-1}}\sum_{q\geq 1}{1\over
q(q+n)}=\sum_{q\geq 1}\frac1q\sum_{n\geq 1}{(-1)^{n-1}\over n^{2a-1}(q+n)}\cdot\eqno(3)$$

We have the following decomposition in the variable $n$:$${1\over n^{2a-1}(q+n)}=\sum_{j=1}^{2a-2}{(-1)^j\over n^{j+1}q^{2a-j-1}}+{1\over
q^{2a-1}}\Big(\frac1n-\frac1{n+q}\Big)\virg\eqno (4)$$
which gives:$$\eqalign{\widetilde\z^*(2a,1)=&\sum_{q\geq 1}\sum_{n\geq 1}(-1)^{n-1}\sum_{j=1}^{2a-2}{(-1)^j\over
n^{j+1}q^{2a-j}}+\sum_{q\geq 1,n\geq 1}{(-1)^{n-1}\over
q^{2a}}\Big(\frac1n-\frac1{n+q}\Big)\cr=&\sum_{j=1}^{2a-2}(-1)^j\sum_{q\geq 1,n\geq 1}{(-1)^{n-1}\over
n^{j+1}q^{2a-j}}+\sum_{q\geq 1,n\geq 1}{(-1)^{n-1}\over
q^{2a}}\Big(\frac1n-\frac1{n+q}\Big)\cr=&\sum_{j=1}^{2a-2}(-1)^j\widetilde\z(j+1)\z(2a-j)+u,}\eqno(5)$$ by defining, for $s>1$,
 $\ds\widetilde\z(s)=\sum_{n\geq 1}\frac{(-1)^{n-1}}{n^s}$  and
$\ds u:=\sum_{q\geq 1,n\geq 1}{(-1)^{n-1}\over
q^{2a}}\Big(\frac1n-\frac1{n+q}\Big)\cdot$ Remark that if $a=1$, the first sum in the last line  of (5) does not exist, as it was
already the case in (4). In the sequel, all sums indexed by a void set are identically 0.  Now:
$$u=\sum_{n\geq 1}\sum_{q\geq 1}{(-1)^{n-1}\over q^{2a-1}n(n+q)}=\sum_{m\geq 1}{(-1)^m\over m}\sum_{q=1}^{m-1}{(-1)^{q-1}\over
q^{2a-1}(m-q)}\eqno (6)$$ by defining $n+q=m$ in the first summation. By new partial fractions decomposition, one has:$${1\over
(m-q)q^{2a-1}}=\sum_{j=1}^{2a-2}{1\over q^{j+1}m^{2a-j-1}}+{1\over m^{2a-1}}\Big(\frac1q+\frac1{m-q}\Big)\virg\eqno(7)$$ which gives
:$$u=\sum_{j=1}^{2a-2}\sum_{m>q\geq 1}{(-1)^{m+q-1}\over q^{j+1}m^{2a-j}}+\sum_{m>q\geq 1}{(-1)^{m+q-1}\over
m^{2a}}\Big(\frac1q+\frac1{m-q}\Big)\cdot\eqno (8)$$ Let  $u'$ be the first sum  of (8) an $u''$ the second one. The sum
$u'$ splits in two
 parts : we sum first  from  1 to $a-1$, and from $a$ to
$2a-2$. In the second part we put $j=2a-1-j'$. Then $j+1=2a-j'$ and $2a-j=j'+1$ : so this second part  becomes  (we  replace index $j'$ by $j$)
$\ds
\sum_{j=1}^{a-1}\sum_{m>q\geq 1}{(-1)^{m+q-1}\over m^{j+1}q^{2a-j}}\virg$ which may be written, by exchange of mute indexes 
$q$ and
$m$, $\ds
\sum_{j=1}^{a-1}\sum_{q>m\geq 1}{(-1)^{m+q-1}\over q^{j+1}m^{2a-j}}\cdot$ So$$u'=\sum_{j=1}^{a-1}\sum_{q\not=m}{(-1)^{m+q-1}\over
q^{j+1}m^{2a-j}\cdot}\cdot\eqno(9)$$

 Now we observe that $$\sum_{q\not=m}{(-1)^{m+q}\over
q^{j+1}m^{2a-j}}=\widetilde\z(j+1)\widetilde\z(2a-j)-\z(2a+1).$$ So  equality (9) becomes
$$u'=-\sum_{j=1}^{a-1}\widetilde\z(j+1)\widetilde\z(2a-j)+(a-1)\z(2a+1).\eqno(10)$$  In fact, we may also write
:$$u'=-\sum_{j=1}^{a-1}\widetilde\z(2j)\widetilde\z(2a+1-2j)+(a-1)\z(2a+1).\eqno(10')$$ 

Now look at the sum $u''$. We have :$$u''=\sum_{m\geq 1}{(-1)^m\over m^{2a}}\mu(m)$$ where
$$\mu(m):=\sum_{q=1}^{m-1}(-1)^{q-1}\Big(\frac1q+\frac1{m-q}\Big)\cdot$$
   Clearly, $\mu(m)$ vanishes if $m$ is odd and
$$\mu(2m)=2\sum_{q=1}^{2m-1}{(-1)^{q-1}\over q}=2(H_{2m-1}-H_{m-1})=2(H_{2m}-H_m)+\frac1m\cdot$$ As a result:
$$\eqalign{u''=&\sum_{m\geq 1}{1\over (2m)^{2a}}\Big(2(H_{2m}-H_m)+\frac1m\Big)\cr=&{1\over 2^{2a-1}}Z(2a)-{1\over
2^{2a-1}}\zeta^*(2a,1)+{1\over 2^{2a}}\z(2a+1),}\eqno (11)$$

Remark that $$\zeta^*(2a,1)-\widetilde\z^*(2a,1)=\sum_{m\geq 1}{1-(-1)^{m-1}\over m^{2a}}H_m={1\over 2^{2a-1}}Z(2a).\eqno (12)$$

So collecting previous results we obtain :$$\eqalign{\widetilde\z^*(2a,1)=
\sum_{j=1}^{2a-2}(-1)^j\widetilde\z(j+1)\z(2a-j)-\sum_{j=1}^{a-1}\widetilde\z(j+1)\widetilde\z(2a-j)+&(a-1)\z(2a+1)\cr+(\zeta^*(2a,1)-\widetilde\z^*(2a,1))-{1\over
2^{2a-1}}\zeta^*(2a,1)+{1\over 2^{2a}}\z(2a+1).}$$
This gives
$$\eqalign{2\widetilde\z^*(2a,1)&=\sum_{j=1}^{2a-2}(-1)^j\widetilde\z(j+1)\z(2a-j)-\sum_{j=1}^{a-1}\widetilde\z(j+1)\widetilde\z(2a-j)\cr&+(1-2^{1-2a})\zeta^*(2a,1)+(a-1+2^{-2a})\z(2a+1).}\eqno(13)$$
At this stage, using (1.4), we can write $$2\widetilde\z^*(2a,1)=A+B,\eqno(14)$$ putting $$\ds
A=\Big(2a-{2a+1\over2^{2a}}\Big)\z(2a+1)=(2a+1)\widetilde\z(2a+1)-\z(2a+1)\eqno(15)$$ 
and
$$B=\sum_{j=1}^{2a-2}(-1)^j\widetilde\z(j+1)\z(2a-j)-\sum_{j=1}^{a-1}\widetilde\z(j+1)\widetilde\z(2a-j)-\Big(1-\frac1{2^{2a-1}}\Big)\sum_{j=1}^{a-1}\z(j+1)\z(2a-j)\eqno
(16)$$ We will give a   reduction of $B$. 
We write first (by consideration of parity) :$$\eqalign{B=&\sum_{j=1}^{a-1}\widetilde\z(2j+1)\z(2a-2j)-\sum_{j=1}^{a-1}\widetilde\z(2j)\z(2a+1-2j)
\cr&-\sum_{j=1}^{a-1}\widetilde\z(2j+1)\widetilde\z(2a-2j)-\Big(1-\frac1{2^{2a-1}}\Big)\sum_{j=1}^{a-1}\z(2j+1)\z(2a-2j)\cr=&\sum_{j=1}^{a-1}{1\over
2^{2a-2j-1}}\widetilde\z(2j+1)\z(2a-2j)-\sum_{j=1}^{a-1}\widetilde\z(2j)\z(2a+1-2j)\cr&-\Big(1-\frac1{2^{2a-1}}\Big)\sum_{j=1}^{a-1}\z(2j+1)\z(2a-2j).}\eqno
(17)$$
By replacing $\widetilde\z(2j+1)$ by $(1-2^{-2j})\z(2j+1)$ in the first sum we obtain,  by grouping with the third one:
$$\eqalign{B=&\sum_{j=1}^{a-1}\Big({1\over
2^{2a-2j-1}}-1\Big)\z(2j+1)\z(2a-2j)-\sum_{j=1}^{a-1}\widetilde\z(2j)\z(2a+1-2j)\cr=&-\sum_{j=1}^{a-1}\widetilde\z(2j+1)\z(2a-2j)-\sum_{j=1}^{a-1}\widetilde\z(2j)\z(2a+1-2j)
\cr=&-2\sum_{j=1}^{a-1}\widetilde\z(2j)\z(2a+1-2j),}\eqno(18)$$ the last equality resulting for the exchange $j\mapsto a-j.$
Using (14) and (15) we get finally:$$\tableau{\widetilde\z^*(2a,1)=\Big(a+\frac12\Big)\widetilde\z(2a+1)-\frac12\z(2a+1)-\sum_{j=1}^{a-1}\widetilde\z(2j)\z(2a+1-2j)\crh}\,.\eqno (19)$$

\noindent {\bf Examples }

$a=1$: $\ds\sum_{n\geq 1}(-1)^{n-1}{H_n\over
n^{2}}=\frac32\widetilde\z(3)-\frac12\z(3)=\Big(\frac32\times\frac34-\frac12\Big)\z(3)=\frac58\z(3) ;$

$a=2$: $\ds\sum_{n\geq 1}(-1)^{n-1}{H_n\over
n^{4}}=\frac52\widetilde\z(5)-\frac12\z(5)-\widetilde\z(2)\z(3)=\frac{59}{32}\z(5)-\frac12\z(2)\z(3).$\vskip 1mm 

In fact  formula (19) was obtained in [2] by Flajolet and Salvy, but they used there  residue's theorem. And it was etablished by Sitaramachandrarao in [5] by
using the result for $J(2a)$ given by Jordan in [4], but we will inverse the processes of Sitaramachandrarao, because the proof in [4] is hard to understand.
\vskip 3mm
\noindent{\bf II.2 The sum $Z(2a)$}\vskip 2mm

 Subtracting $(19)$ from $(1.4)$ and using (12) we get $$2^{-2a+1}Z(2a)=$$
$$(a+1)\z(2a+1)-\sum_{j=1}^{a-1}\z(2j)\z(2a+1-2j)-\Big(a+\frac12\Big)\widetilde\z(2a+1)+\frac12\z(2a+1)+\sum_{j=1}^{a-1}\widetilde\z(2j)\z(2a+1-2j)$$
which gives $$2^{-2a+1}Z(2a)
=\Big(a+\frac32\Big)\z(2a+1)-\Big(a+\frac12\Big)\widetilde\z(2a+1)-\sum_{j=1}^{a-1}{1\over
2^{2j-1}}\z(2j)\z(2a+1-2j).$$Replacing
$\widetilde\z(2a+1)$ by its expression in terms of $\z(2a+1)$, one obtains now :$$\tableau{\sum_{n\geq 1}{H_{2n}\over
n^{2a}}=\frac14(2a+1+2^{2a+1})\z(2a+1)-\sum_{j=1}^{a-1}2^{2a-2j}\z(2j)\z(2a+1-2j)\crh}\,.\eqno (20)$$ This last expression is  given in [5], Theorem  1, but
using Bernoulli polynomials and integral transformations. This relation (20) will be the cornerstone for our proof of first Jordan's formula.

Relations (20) et $(1.4)$ give easily the formula:$$\tableau{\ds \sum_{n\geq 1}{H_{2n-1}\over
(2n-1)^{2a}}={2a+1\over 2}\lambda(2a+1)-\sum_{j=1}^{a-1}\z(2a+1-2j)\lambda(2j)\crh}\,\,,\eqno(20')$$ 
\vskip 3mm
\noindent{\bf III. Proof of the first Jordan formula}\vskip 2mm
In this section, we  give a close formula for the sum $J(2a)$ in terms of $\z$ and $\lambda $ series, first by using previous results and second by an integral
representation of it.\vskip 2mm

\noindent{\bf III.1
First proof by sums of series }\vskip 3mm

Clearly, we have $$H_{2n}=S_n+\frac12 H_n \,,$$which gives $$\sum_{n\geq 1}{H_{2n}\over
n^{2a}}=J(2a)+\frac12 \zeta^*(2a,1).\eqno (21)$$
By
(20),(21), and $(1.4)$ we obtain
:$$\eqalign{J(2a)=&\frac14(2a+1+2^{2a+1})\z(2a+1)-\sum_{j=1}^{a-1}2^{2a-2j}\z(2j)\z(2a+1-2j)\cr&-\frac12((a+1)\z(2a+1))+\frac12\sum_{j=1}^{a-1}\z(2j)\z(2a+1-2j)\cr=&{2^{2a+1}-1\over
4}\z(2a+1)-\frac12\sum_{j=1}^{a-1}(2^{2a+1-2j}-1)\z(2j)\z(2a+1-2j).}$$ The transformation $j\mapsto a-j$ furnishes the final expression:$$\tableau{J(2a)=\sum_{n\geq 1}{S_n\over
n^{2a}}={2^{2a+1}-1\over 4}\z(2a+1)-\frac12\sum_{j=1}^{a-1}(2^{2j+1}-1)\z(2j+1)\z(2a-2j)\crh}\,.\eqno(22)$$
This formula becomes, by using $\lambda$'s series: 
$$\tableau{J(2a)=\sum_{n\geq 1}{S_n\over
n^{2a}}=2^{2a-1}\lambda(2a+1)-\sum_{j=1}^{a-1}2^{2j}\lambda(2j+1)\z(2a-2j)\crh}\,.\eqno(23)$$

\noindent{\bf Examples}

 $\ds \sum_{n\geq 1}{S_n\over n^{2}}=\frac74\z(3)\,\,\,\,;\,\,\,\sum_{n\geq 1}{S_n\over n^4}=\frac{31}4\z(5)-\frac72\z(3)\z(2).$
\vskip 2mm

 \noindent{\bf III.2 Second method for calculating   $J(2a)$, by an integral representation.}\vskip 2mm

We base this new method on the following claim:
  \proclaim Claim. For $n\geq 0,$ we have the relation  $\displaystyle 2S_n=\int_0^\pi{\sin^2nx\over \sin x}\,\,{\rm d}x.$\par {\bf Proof.} Denote by $I_n$ the integral in
the claim. For $n\geq 1$,
 by  the classical relation 
$\sin^2nx-\sin^2(n-1)x=\sin x\sin(2n-1)x,  $ we have:$$I_n-I_{n-1}=\int_0^\pi \sin(2n-1)x\, {\rm d}x=\frac2{2n-1}\virg$$ and this prove our claim, since
$I_0=0$.\cqfd  \vskip 1mm

 It follows from the claim the crucial integral representation for $J(2a)$:
$$2J(2a)=\int_0^\pi{\varphi_a(x)\over
\sin x}\,\,{\rm d}x,$$ where$$\varphi_a(x):=\sum_{p\geq 1}{\sin^2px \over p^{2a}}\cdot\eqno (24) $$ The function $\varphi_a$ is expressible by    Bernoulli
polynomials, but we did not use them explicitly. In fact, we   use  successive     partial integrations in $(24)$, which correspond to  induction relations between
this polynomials. Put $$\a_n(a)=\int_0^\pi{\varphi_a(x)\over \sin x}\cos 2nx\, {\rm d}x .$$
We have $2J(2a)=\a_0(a)$. For all $a\geq 1$, the functions $\ds\psi_a: x\mapsto {\varphi_a(x)\over \sin x}$ are integrable on the interval $[0,\pi]$. For $a\geq 2$,
this results from the inequality $\abs{\sin px}\leq p\abs{\sin x},$ which gives,  $$\abs {\psi_a(x)}\leq \z(2a-1)<+\infty.$$ And for $a=1$, the integrability results
from  the classical expansion (for $0\leq x\leq \pi)$:$$\frac12x(\pi-x)=\sum_{p\geq 1}{\sin^2px\over p^2}=\varphi_2(x)$$ and from the fact that $\psi_2(x)  $ tends to $\ds \frac\pi2$ when
$x$ tends to 0 or to
$\pi$.
Hence, by  Riemann-Lebesgue theorem, the sequence      
$\a_n(a)$ tends to 0 if $n$ goes to infinity. 

 For $n\geq 1,$ we have :$$\a_{n-1}(a)-\a_n(a)=\int_0^\pi{\varphi_a(x)\over \sin x}(\cos(2n-2)x-\cos2nx){\rm
d}x=2\int_0^\pi\varphi_a(x)\sin(2n-1)x\,\,{\rm d}x:=2v_n.$$ So, by summation, $\ds \a_0(a)-\a_n(a)=2\sum_{k=1}^nv_k$. Then, $\displaystyle \a_0(a)=2\sum_{n=1}^\infty v_n,$ and finally 
$$J(a)=\sum_{n=1}^\infty v_n .$$

In the following, we drop the subscript $a$ in $\varphi_a.$ \vskip 2mm

\noindent{\bf Calculation of} $\ds v_n=\int_0^\pi\varphi (x)\sin(2n-1)x\,\,{\rm
d}x$. \vskip 1mm

We proceed by successive integrations by part, involving the derivatives of $\varphi$
:$$\varphi^{(2j)}(x)=(-1)^{j-1}2^{2j-1}\sum_{k=1}^\infty{\cos 2kx\over k^{2a-2j}}\virg\,\,\,\, \varphi^{(2j-1)}(x)=(-1)^{j-1}2^{2j-2}\sum_{k=1}^\infty{\sin 2kx\over k^{2a-2j+1}}$$

when  $1\leq j\leq a-1$ for  even derivatives and $1\leq j\leq a$ for  odd ones. We observe that if $j=a,$ the second relation makes sense  only if 
$x\not=0
\,\,\,({\rm mod}\,\,\,\pi)
$. In this case we have
$$\varphi^{(2a-1)}(x)=(-1)^{a-1}2^{2a-2}\sum_{k=1}^\infty{\sin 2kx\over k}\cdot$$ Using Fourier's series or complex logarithm, this is (when
$0<x<\pi$):
$$\varphi^{(2a-1)}(x)=(-1)^{a-1}2^{2a-3}(\pi-2x).$$ 

In the following calculation, $n$ appears only by mean of $2n-1$, so   we set  $2n-1=m$.
Since 
$\varphi(0)=\varphi(\pi)=0$, we obtain first :
$$v_n=\Big[-{\varphi(x)\cos mx\over m }\Big]_0^\pi+\frac1m\int_0^\pi\varphi'(x)\cos mx \,{\rm d}x=\frac1m\int_0^\pi\varphi'(x)\cos mx \,{\rm d}x.$$  When $a=1$, this  process stops
$\ds(\varphi'(x)=\frac12(\pi-2x))$ and we obtain :$$v_n=\frac1{m}\int_0^\pi \Big(\frac\pi2-x\Big)\cos mx \,{\rm d}x=\frac1{m^2}\Big[ \Big(\frac\pi2-x\Big)\sin mx \Big]_0^\pi
+\frac1{m^2}\int_0^\pi
\sin mx\,{\rm d}x=\frac2{m^3}\eqno(24')$$ For $a\geq 2,$   we use two   integrations by part :$$\eqalign {v_n&=\frac1{m^2}\Big[\varphi'(x)\sin mx\Big]_0^\pi -\frac1{m^2}\int_0^\pi
\varphi''(x)\sin mx\,{\rm d}x=-\frac1{m^2}\int_0^\pi\varphi''(x)\sin mx\,{\rm d}x\cr &={1\over m^3}\Big[\varphi''(x)\cos mx\Big]_0^\pi-{1\over
m^3}\int_0^\pi \varphi'''(x)\cos mx\,{\rm d}x\cr &=-{4\over m^3}\z(2a-2)-{1\over
m^3}\int_0^\pi \varphi'''(x)\cos mx\,{\rm d}x,}$$ the last equality coming from  the expression of $\varphi''$ (remember that  $m$ is odd). When 
$a=2$, we stop. Eventually, we stop when the last integral is
$\ds\int_0^\pi\varphi^{(2a-1)}(x)\cos mx {\rm d}x$ and obtain the relation:$$v_n=-\sum_{j=1}^{a-1}{2^{2j}\over m^{2j+1}}\z(2a-2j)+{(-1)^{a-1}\over
m^{2a-1}}\int_0^\pi\varphi^{(2a-1)}(x)\cos mx {\rm d}x.$$ 
 which gives, by
taking the value of 
 $\varphi^{(2a-1)}(x)$ into account and the previous relation $(24')$:
$$v_n=-\sum_{j=1}^{a-1}{2^{2j}\over m^{2j+1}}\z(2a-2j)+{2^{2a-1}\over m^{2a+1}}\cdot$$

Finally, we obtain :  $$J(2a)=2^{2a-1}\lambda(2a+1)-\sum_{j=1}^{a-1}2^{2j}\lambda(2j+1)\z(2a-2j),$$
which is the already seen formula $(23)$.\vskip 3mm

\noindent{\bf IV. Calculation of $\ds \bar J(2a)=\sum_{n\geq 1}{S_n\over(2 n-1)^{2a}}$ and applications}\vskip 4mm
\noindent{\bf IV.1 Calculation of $\ds \bar J(2a)$}\vskip 2mm
First we write a relation between $\bar J(b)$ and $J(b)$ for all integer $b\geq 2.$
We start by the following easy identity :$$S_n=\sum_{p\geq 1}{2n\over(2p-1)(2n+2p-1)}=\sum_{p\geq 1}\Big({1\over 2p-1}-{1\over 2n+2p-1}\Big)\cdot\eqno(25)$$
 So, we have$$\bar J(b)=\sum_{n\geq 1,p\geq 1}{2n\over
(2n-1)^b(2p-1)(2n+2p-1)}=\sum_{p\geq 1}\frac1{2p-1}\sum_{n\geq 1}{2n\over
(2n-1)^b(2n+2p-1)}\cdot\eqno(26)$$  
We denote by  $\rho(p) $ the inner summation  in (26) and put $\ds Q(x,\mu)={x+1\over x^b(x-\mu)}\cdot $ So the  generic term in   $\rho (p)$ is
$Q(2n-1,-2p).$ We decompose the rational fraction $Q(x,\mu)$ into  partial fractions in the variable   $x$:
$$Q(x,\mu)=-{1\over \mu x^b}-\sum_{j=1}^{b-2}{\mu+1\over\mu^{j+1}x^{b-j}}-{\mu+1\over \mu^b}\Big(\frac1x-\frac1{x-\mu}\Big),$$ so that
$$\rho(p)=\sum_{n\geq 1}\Big({1\over 2p(2n-1)^b}+\sum_{j=1}^{b-2}{(-1)^{j+1}(2p-1)\over (2p)^{j+1}(2n-1)^{b-j} }+(-1)^b{2p-1\over
(2p)^b}\Big(\frac1{2n-1}-\frac1{2n+2p-1}\Big)\Big)\eqno(27)$$
or  using (25) $$\rho(p)=\frac1{2p}\lambda(b)+(2p-1)\sum_{j=1}^{b-2}{(-1)^{j+1}\over (2p)^{j+1} }\lambda(b-j)+(-1)^b{2p-1\over
(2p)^b}S_p.$$By carrying this value in (26), we get:$$\bar
J(b)=\lambda(b)\sum_{p\geq 1}\frac1{2p(2p-1)}+\sum_{j=1}^{b-2}(-1)^{j+1}\lambda(b-j)\sum_{p\geq 1}\frac1{(2p)^{j+1}}+(-1)^b{S_p\over (2p)^b}\cdot$$
The first sum in the previous line is  $\ln 2$. One gets the general relation, for all integer $b$:$$\tableau{\bar
J(b)+{(-1)^{b-1}\over 2^b}J(b)=\lambda (b)\ln 2+\sum_{j=1}^{b-2}{(-1)^{j+1}\over 2^{j+1}}\lambda(b-j)\z(j+1)\crh}.\eqno(28)$$

So we can calculate $\bar J(b)$ if   $J(b)$ is known, in particulier when $b$ is even or  when $b=3$ (see remark below). Look at $b=2a$. We use
then
$(23)$ by replacing 
$j$ by $a-j$ in the sum:
$$\eqalign{\bar J(2a)=&\frac12\Big(1-{1\over 2^{2a+1}}\Big)\z(2a+1)+\lambda (2a)\ln2-{1\over
2^{2a+1}}\sum_{j=1}^{a-1}2^{2a-2j+1}\lambda(2a-2j+1)\z(2j)\cr&+\sum_{j=1}^{2a-2}{(-1)^{j+1}\over 2^{j+1}}\lambda(2a-j)\z(j+1).}\eqno(29)$$ We split the second
sum in $(29)$ accordingly with the parity of  the index $j$.  This gives  :$$\eqalign{\sum_{n\geq 1}{S_n\over
(2n-1)^{2a}}=&\lambda(2a)\ln2+\frac12\Big(1-{1\over 2^{2a+1}}\Big)\z(2a+1)-\sum_{j=1}^{a-1}{1\over
2^{2j+1}}\lambda(2a-2j)\z(2j+1)\cr},$$ and finally: $$\tableau{\bar J(2a)=\sum_{n\geq 1}{S_n\over
(2n-1)^{2a}}=\lambda(2a)\ln2+\frac12\lambda(2a+1)-\sum_{j=1}^{a-1}{1\over 2^{2j+1}}\lambda(2a-2j)\z(2j+1)\crh}\,\,.\eqno (30)=(b)$$
This formula
is given  by Jordan in $[4]$ and is used in [5].

\noindent{\bf Examples :}\vskip 2mm

$\ds \bar J(2)=\frac34\z(2)\ln2+\frac7{16}\z(3),\,\,\,\hskip 5mm \bar J(4)=\frac{15}{16}\z(4)\ln2+{31\over 64}\z(5)-{3\over 32}\z(2)\z(3)$\vskip 2mm

\noindent{\bf The case $b=3$}  

By using the classical series $\ds {\rm li}_4\frac12=\sum_{n\geq 1}{1\over 2^nn^4}$ the sums $\widetilde \zeta^*(3,1)$ (and then $J(3)$) can be evaluated (see for instance
[2]) :$$\widetilde \zeta^*(3,1)=-2{\rm li}_4\frac12+\frac{11}4\z(4)+\frac12\z(2)\ln^22-\frac1{12}\ln^42-\frac74\z(3)\ln 2.$$ This gives
first:$$J(3)=\sum_{n\geq 1}{H_{2n}\over n^3}-\frac12\zeta^*(3,1)=4(\zeta^*(3,1)-\widetilde \zeta^*(3,1))-\frac12\zeta^*(3,1)=\frac72\zeta^*(3,1)-4\widetilde \zeta^*(3,1)    .$$ But we
know
$\ds
\zeta^*(3,1)=\frac54\z(4)$ (formula (1.3) for $b=3$). So we obtain $$J(3)=8{\rm li}_4\frac12-\frac{53}8\z(4)-2\z(2)\ln^22+\frac13\ln^42+7\z(3)\ln 2.$$
Then, formula (28) gives, when $b=3$:$$\bar J(3)+\frac18J(3)=\lambda(3)\ln2+\frac14\lambda(2)\z(2)=\frac78\z(3)\ln 2+\frac{15}{32}\z(4), $$ this gives:$$\bar
J(3)=-{\rm li}_4\frac12+\frac{83}{64}\z(4)+\frac14\z(2)\ln^22-\frac1{24}\ln^42.$$
 Infortunatly,  it is conjectured that 
$\widetilde \zeta^*(k,1)$ when $k\geq 5$ is odd are not expressible by usual constants and so probably the $J(b)$ are not simple when $b$ is odd $\geq 5$.
\vskip 3mm

\noindent{\bf III.2 Application to the calculation of new sums involving the $H_n$'s}
\vskip 3mm
We establish  a pretty formula:  by using  $(20')$ (in which we change  the $j$ index by $a-j$), previous relation  
 (30) and the simple fact that $$H_{2n-1}=\frac12H_{n-1}+S_n$$ we obtain$$\frac12\sum_{n\geq 1}{H_{n-1}\over(2n-1
)^{2a}}=-\lambda(2a)\ln2+a\lambda(2a+1)+\sum_{j=1}^{a-1}\lambda(2j)\zeta(2a+1-2j)\Big({1\over2^{2a+1-2j}}-1\Big),$$ or
$$\tableau{\sum_{n\geq 1}{H_{n}\over (2n+1
)^{2a}}=-2\lambda(2a)\ln2+2a\lambda(2a+1)-2\sum_{j=1}^{a-1}\lambda(2j)\lambda(2a+1-2j)\crh}\,.\eqno (31)=(c)$$\vskip 2mm Thus, the closed form of the previous sum
  involves only  
$\lambda$ series.
\vskip 2mm
\noindent{\bf Examples:}\vskip 2mm 

$\ds \sum_{n\geq 1}{H_{n}\over (2n+1
)^2}=-{\pi^2\over 4}\ln 2+2\lambda(3),\hskip 5mm\sum_{n\geq 1}{H_{n}\over (2n+1
)^4}=-{\pi^4\over 48}\ln 2+4\lambda(5)-\frac{\pi^2}4\lambda(3).$

\vskip 4mm

\noindent{\bf IV. Series including the sums }$\ds \sum_1^n{(-1)^{k-1}\over k^{2a}}$

\vskip 4mm
Now we define, for  $b\geq 1$ $\widetilde H_0^{(b)}=0$ and, if $n\geq 1$  $$\widetilde H_{n}^{(b)}=\sum_{k=1}^n {(-1)^{k-1}\over k^b}\cdot$$
We wish  to get a formula for the sum$$\sum_{n=1}^\infty
{(-1)^n\widetilde H_{n-1}^{(2a)}\over n}\eqno (32)$$
when $a$ is a positive integer. This sum exists because we can write $\widetilde H_{n-1}^{(2a)}=\widetilde \z(2a)+r_n$, where the ``remainder" $r_n={\rm
O}({1\over n^{2a}})$ and the series (32) is the sum of two convergent series.

 Return to the second line of (5) and its second sum.
We have, by absolute convergence $$\sum_{n\geq 1,q\geq 1}{(-1)^{n-1}\over
q^{2a}}\Big(\frac1n-\frac1{n+q}\Big)=\lim_{r\tend \infty}\Bigg(\sum_{k=1}^r\sum_{q+n=k}{(-1)^{n-1}\over
q^{2a}n}-\sum_{k=1}^r\sum_{ q+n=k}{(-1)^{n-1}\over
q^{2a}(n+q)}\Bigg).\eqno(33)$$ But the first term in (33) is the 
Cauchy product of two convergent series, one of them being absolutly convergent and so the limit
for $r$ infinite is $\ln 2\z(2a)$. Thus we have :$$\widetilde\z^*(2a,1)=\ln2\z(2a)+\sum_{j=1}^{2a-2}(-1)^j\widetilde\z(j+1)\z(2a-j)-\sum_{k=1}^\infty\sum_{
q+n=k}{(-1)^{n-1}\over q^{2a}(n+q)}\cdot$$

But $\ds \sum_{ q+n=k}{(-1)^{n-1}\over
q^{2a}(n+q)}=-{(-1)^{k-1}\over k}\sum_{q=1}^{k-1}{(-1)^{q-1}\over q^{2a}}=-{(-1)^{k-1}\widetilde H_{k-1}^{(2a)}\over k}\virg$
which gives $$\widetilde\z^*(2a,1)=\ln2\z(2a)+\sum_{j=1}^{2a-2}(-1)^j\widetilde\z(j+1)\z(2a-j)+\sum_{k=1}^\infty {(-1)^{k-1}\widetilde H_{k-1}^{(2a)}\over k}
\cdot$$ By the expression $(19)$ of   $\widetilde\z^*(2a,1)$  we get:$$\eqalign{\sum_{n=1}^\infty {(-1)^n\widetilde H_{n-1}^{(2a)}\over
n}=&\ln2\z(2a)+\sum_{j=1}^{2a-2}(-1)^j\widetilde\z(j+1)\z(2a-j)\cr&-\Big(a+\frac12\Big)\widetilde\z(2a+1)+\frac12\z(2a+1)+\sum_{j=1}^{a-1}\widetilde\z(2j)\z(2a+1-2j).}$$

Separating odd indices $j$  from even one's in the first sum of right hand member, we obtain:
$$\tableau{\eqalign{\sum_{n=1}^\infty {(-1)^n\widetilde H_{n-1}^{(2a)}\over
n}&=\frac12\z(2a+1)-\big(a+\frac12\big)\widetilde\z(2a+1)\cr&+\ln2\z(2a)+\sum_{j=1}^{a-1}\widetilde\z(2j+1)\z(2a-2j)}\crh }\,.\eqno(34)$$

\noindent {\bf Examples}

$\ds \sum_{n=1}^\infty
{(-1)^n\widetilde H_{n-1}^{(2)}\over n}=\ln 2\z(2)+\frac12\z(3)-\frac32\widetilde\z(3)=\ln 2\z(2)-\frac58\z(3).$

$\ds \ds \sum_{n=1}^\infty
{(-1)^n\widetilde H_{n-1}^{(4)}\over n}=\ln 2\z(4)-\frac{59}{32}\z(5)-\frac34\z(3)\z(2).$

\vskip 6mm

\noindent{\bf V. Other new sums}

We define for all integers $t\geq 1$ and $s\geq 2$, $\ds S_n^{(t)}=\sum_{k=1}^n\frac1{(2k-1)^t}$ and \vskip 2mm

 $$\ds \sigma(s,t)=\sum_{n\geq 1}{S_n^{(t)}\over n^s} $$ (whose
``weight" is $s+t$) and we investigate some relations between these new sums, by imitation of well known one's regarding the $\ds
\z(s,t)=\sum_{n\geq 1}{H_{n-1}^{(t)}\over n^s}\cdot$  For these last sums,  results was first given by Euler for odd weights
$\leq 13$ and more recently  established for all odd weights by D. Borwein, J.M. Borwein and R. Girgensohn in [1] and Flajolet, Salvy in [2].\vskip 2mm

\noindent{\bf V.1 Linear relations between  $\sigma (s,t)$ of same weight $s+t$ and $\ds\sum_{n\geq 1}  {H_{n}\over (2n+1
)^{s+t-1}}$ }

  We write  $$\sigma(s,t)=\sum_{n\geq m\geq 1}{1\over
n^s(2m-1)^t}=\sum_{p\geq 0}\sum_{m\geq 1}{1\over (m+p)^s(2m-1)^t}\cdot$$ Recall the decomposition in ${\bb C}[[u]]$:$${1\over (r-u)^tu^s}=\sum_{i=0}^{s-1}{t+i-1\choose
i}{1\over r^{t+i}u^{s-i}}+\sum_{j=0}^{t-1}{s+j-1\choose j}{1\over r^{s+j}(r-u)^{t-j}}\virg\eqno(35)$$ where $r$ is some parameter. We set $\ds u=m+p, r=\frac12+p$
and we  get, by (35) :$$\eqalign{{1\over
(m+p)^s(2m-1)^t}=&(-1)^t\sum_{i=0}^{s-1}{t+i-1\choose i}{2^{i}\over (2p+1)^{t+i}(m+p)^{s-i}}\cr&+2^s\sum_{j=0}^{t-1}{s+j-1\choose
j}{(-1)^j\over(2p+1)^{s+j}(2m-1)^{t-j}}\virg}$$We cut the right hand of previous equality in three parts :$$A(m,p)=(-1)^t\sum_{i=0}^{s-2}{t+i-1\choose
i}{2^{i}\over (2p+1)^{t+i}(m+p)^{s-i}}\virg$$ $$B(m,p)=2^s\sum_{j=0}^{t-2}{s+j-1\choose
j}{(-1)^j\over(2p+1)^{s+j}(2m-1)^{t-j}}$$ and the part corresponding to $i=s-1$ and $j=t-1$ (then the corresponding binomial coefficients are equal):
$$C(m,p)=(-1)^t2^{s-1}{s+t-2\choose
s-1}\times{1\over(2p+1)^{s+t-1}}\Big({1\over m+p}-{2\over 2m-1}\Big)$$(note that the sum of residue must  be 0 ). Now we
have$$\sigma(s,t)=\sum_{p\geq 0,m\geq 1}A(m,p)+B(m,p)+C(m,p),$$$$\sum_{p\geq 0,m\geq 1}A(m,p)=(-1)^t\sum_{i=0}^{s-2}2^{i}{t+i-1\choose i}\sigma(s-i,t+i),$$
$$\sum_{p\geq 0,m\geq 1}B(m,p)=2^s\sum_{j=0}^{t-2}(-1)^j{s+j-1\choose
j}\lambda(s+j)\lambda(t-j).$$ In order to evaluate the last term, we write $${1\over m+p}-{2\over 2m-1}={1\over m+p}-\frac1m+2\Big({1\over 2m}-{1\over
2m-1}\Big).$$ By summing (in $m$) we get $-H_p-2\ln 2$, for any $p\geq 0$ and so $$\eqalign{\sum_{p\geq 0,m\geq 1}C(m,p)=&(-1)^{t-1}2^{s-1}{s+t-2\choose
s-1}\sum_{p\geq 1}{H_p\over(2p+1)^{s+t-1}}\cr&+(-1)^{t-1}2^s{s+t-2\choose
s-1}\lambda(s+t-1)\,\ln 2.}$$It follows the general relation$$\tableau{\eqalign{(-1)^t\sigma(s,t)=&\sum_{i=0}^{s-2}2^{i}{t+i-1\choose
i}\sigma(s-i,t+i)\cr&+(-1)^t2^s\sum_{j=0}^{t-2}(-1)^j{s+j-1\choose j}\lambda(s+j)\lambda(t-j)\cr&-2^{s-1}{s+t-2\choose
s-1}\sum_{p\geq 1}{H_p\over(2p+1)^{s+t-1}}-2^s{s+t-2\choose
s-1}\lambda(s+t-1)\ln 2}\crh}\eqno (36)$$

If $t=1$ and $s\geq 2$, the relation (36) reduces to $$\sigma(s,1)+\sum_{i=1}^{s-2}2^{i-1}\sigma(s-i,1+i)=2^{s-2}\sum_{p\geq 0}{H_p\over
(2p+1)^s}+2^{s-1}\lambda(s)\ln2.\eqno(36.1)$$
When  $s=2$, it gives $\ds \sigma(2,1)=\sum_{p\geq 1}{H_p\over
(2p+1)^2}+2\lambda(2)\ln2$, which   results also from (23) and (32.1).
\vskip 1mm

When $t$ is even, the left hand member and the first term of right hand member cancel out
and we can write, for all  $t=2r, r\geq 1$:$$\tableau{\eqalign{&\sum_{i=1}^{s-2}2^{i-1}{2r+i-1\choose
i}\sigma(s-i,2r+i)+2^{s-1}\sum_{j=0}^{2r-2}(-1)^j{s+j-1\choose j}\lambda(s+j)\lambda(2r-j)\cr&-2^{s-2}{s+2r-2\choose
s-1}\sum_{p\geq 0}{H_p\over(2p+1)^{s+2r-1}}-2^{s-1}{s+2r-2\choose s-1}\lambda(s+2r-1)\,\ln 2=0}\crh}\eqno (37)$$
{\bf V.2 Application to sums $\ds\sum_{p\geq 1} {H_p\over(2p+1)^{2r+1}}$}

   For example, when $s=2$ and $r=1$ we get the relation$$0=4\lambda(2)^2-4\sum_{p\geq 0}{H_p\over(2p+1)^{3}}-8\ln2\lambda(3).$$But $\ds\lambda(2)^2={\pi^4\over 64}$ et
$\ds
\lambda (3)=\frac78\z(3).$ So, we obtain: $$\sum_{p\geq 1}{H_p\over(2p+1)^{3}}=\lambda(2)^2-2\lambda(3)\ln2={\pi^4\over 64}-\frac74\z (3)\ln 2 .\eqno(38)$$More generally,
$s=2$ gives following  formula  (compare with formula (32.1))$$\sum_{p\geq 1}{H_p\over(2p+1)^{2r+1}}=-2\lambda(2r+1)\ln
2+\frac1r\sum_{j=0}^{2r-2}(-1)^j(j+1)\lambda(2+j)\lambda(2r-j).$$ In order to transform the previous sum
 for $r\geq 2$, denote it  by $D$  and set
$k=2r-2-j$. We get$$D=\sum_{k=0}^{2r-2}(-1)^k(2r-1-k)\lambda(2r-k)\lambda(2+k)=\sum_{j=0}^{2r-2}(-1)^j(2r-1-j)\lambda(2r-j)\lambda(2+j)$$ Writing $\ds
D=\frac12(D+D)$, we get

$$\ds \frac
Dr=\sum_{j=0}^{2r-2}(-1)^j\lambda(2+j)\lambda(2r-j)=\sum_{q=1}^r\lambda(2q)\lambda(2r-2q+2)-\sum_{q=1}^{r-1}\lambda(2q+1)\lambda(2r-2q+1).$$But the first sum in the right
hand side is $\ds \Big(r+\frac12\Big)\lambda(2r+2)$. This last fact follows from the
identity (valid for all integer   $n\geq 1$):$$\tableau{\sum_{j=1}^{n-1}\lambda(2j)\lambda(2n-2j)=\Big(n-\frac12\Big)\lambda(2n)\crh}\,,\eqno (E)$$ This classical identity can
be proved by using the expansion (for $\abs x<1)$:$${\pi x\over 4}\tan {\pi x\over 2}=\sum_{n\geq 1}\lambda(2n)x^{2n}.$$  So, by setting
$r=a-1$, we obtain:

$$\tableau{\sum_{p\geq 1}{H_p\over(2p+1)^{2a-1}}=-2\lambda(2a-1)\ln
2+\Big(a-\frac12\Big)\lambda(2a)-\sum_{q=1}^{a-2}\lambda(2q+1)\lambda(2a-2q-1)\crh}\,.\eqno (39)=(d)$$ 
When $a=2b$ the relation  (39) becomes: 
$$\tableau{\sum_{p\geq 1}{H_p\over(2p+1)^{4b-1}}=-2\lambda(4b-1)\ln
2+\Big(2b-\frac12\Big)\lambda(4b)-2\sum_{q=1}^{b-1}\lambda(2q+1)\lambda(4b-2q-1)\crh}\,.\eqno (39.1)=(e)$$ and when $a=2b+1$, it becomes 
$$\tableau{\eqalign{\sum_{p\geq 1}{H_p\over(2p+1)^{4b+1}}=&-2\lambda(4b+1)\ln
2+\Big(2b+\frac12\Big)\lambda(4b+2)\cr&-\lambda^2(2b+1)-2\sum_{q=1}^{b-1}\lambda(2q+1)\lambda(4b-2q+1)}\crh}\,.\eqno(39.2)=(e')$$

\noindent{\bf Examples :} 

(1) If $b=1$(39.1) gives $\ds \sum_{p\geq 1}{H_p\over(2p+1)^3}=-2\lambda(3)\ln2+\frac32\lambda(4)$, which is the same thing as (38).

(2) If $b=2$ (39.1) gives $\ds \sum_{p\geq 1}{H_p\over(2p+1)^7}=-2\lambda(7)\ln2+\frac72\lambda(8)-2\lambda(3)\lambda(5).$

(3) If $b=1$ (39.2) gives $\ds \sum_{p\geq 1}{H_p\over(2p+1)^5}=-2\lambda(5)\ln 2+\frac52\lambda(6)-\lambda^2(3)$

  and for
$b=2$:$$\sum_{p\geq 1}{H_p\over(2p+1)^9}=-2\lambda(9)\ln2+\frac92\lambda(10)-2\lambda(3)\lambda(7).$$

\noindent{\bf V.3 Linear relations between the $\sigma$'s of same weight}
\vskip 2mm

\noindent{\bf (1) A general relation }

(a) In (37) we put $s=2v$ and replace the sum $\ds\sum_{p\geq 1}{H_p\over(2p+1)^{2v+2r-1}}$ in it by the value given by (39) when $a=v+r$.  So we obtain
(remark that the term with $\ln 2$ cancels out)$$\tableau{\eqalign{& \sum_{i=1}^{2v-2}2^{i-1}{2r+i-1\choose
i}\sigma(2v-i,2r+i)\cr&+2^{2v-1}\sum_{j=0}^{2r-2}(-1)^j{2v+j-1\choose j}\lambda(2v+j)\lambda(2r-j)\cr& -2^{2v-3}{2v+2r-2\choose
2v-1}(2v+2r-1)\lambda(2v+2r)\cr&+2^{2v-2}{2v+2r-1\choose 2v-1}\sum_{j=1}^{r+v-2}\lambda(2j+1)\lambda(2r+2v-2j-1)=0}\crh}\eqno(40.1)$$

(b) Now we put $s=2v+1$ in (37), by use of $(32.1)$  (in which we put $a=v+r$) ; we obtain similarily:$$\tableau{\eqalign{&
\sum_{i=1}^{2v-1}2^{i-1}{2r+i-1\choose i}\sigma(2v+1-i,2r+i)\cr&+2^{2v}\sum_{j=0}^{2r-2}(-1)^j{2v+j\choose j}\lambda(2v+1+j)\lambda(2r-j)\cr&
-2^{2v}{2v+2r-1\choose 2v}(v+r)\lambda(2v+2r+1)\cr&+2^{2v}{2v+2r-1\choose 2v
}\sum_{j=1}^{r+v-1}\lambda(2j)\lambda(2r+2v-2j+1)=0}\crh}\eqno(40.2)$$

\noindent{\bf (2) Application to a formula} for $\sigma(2,2a-1)$\vskip 2mm

In (36) we set $s=2$ and $t=2a-1$, where $a$ is an integer $\geq 2$. This furnishes, by setting$$h_q: =\sum_{p\geq 1}{H_p\over
(2p+1)^q}$$

$$\ds \sigma(2,2a-1)=2\sum_{j=0}^{2a-3}(-1)^j(j+1)\lambda(j+2)\lambda(2a-1-j)+(2a-1)[h_{2a}+2\lambda(2a)\ln2].\eqno (41.1)$$ The sum in square brackets is given by
(32.1) and we
obtain:$$\eqalign{\frac12\sigma(2,2a-1)=\sum_{j=0}^{2a-3}(-1)^j(j+1)\lambda(j+2)\lambda(2a-1-j)\cr+a(2a-1)\lambda(2a+1)-(2a-1)\sum_{j=1}^{a-1}\lambda(2j)\lambda(2a+1-2j).}\eqno
(41.2)$$By  splitting the  first sum in (41.2) in  $\ds \sum_{j=1}^{a-1}(2j-1)\lambda(2j)\lambda(2a+1-2j)$ corresponding to even indices
and
$\ds-\sum_{j=1}^{a-1}(2a-2j)\lambda(2a+1-2j)\lambda(2j)$ corresponding to odd indices,
(41.2)becomes (after changing $j$ into $a-j$)$$\tableau{\sigma(2,2a-1)=2a(2a-1)\lambda(2a+1)-8\sum_{j=1}^{a-1}j\lambda(2a-2j)\lambda(2j+1)\crh}\,.\eqno (42)=(f)$$
 \noindent {\bf (3) Examples and applications}

$a=1$ : $\sigma(2,1)=2\lambda(3)=\ds \frac74\z(3)$

$a=2$ : $\ds\sigma (2,3)=12\lambda(5)-8\lambda(2)\lambda(3)=\frac{93}8\z(5)-\frac{21}4\z(2)\z(3)$ \vskip 1mm

$a=3$ : $\sigma(2,5)=30\lambda(7)-8\lambda(4)\lambda(3)-16\lambda(5)\lambda(2)$

If we test  $s=3$ in (36.1), we get $\ds \sigma(3,1)+\sigma(2,2)=2\sum_{p\geq 1}{H_p\over
(2p+1)^3}+4\lambda(3)\ln2,$ that is$$\sigma(3,1)+\sigma(2,2)=3\lambda(4)={\pi^4\over 32}=\frac{45}{16}\z(4).\eqno (43)$$We have already seen
that$$J(3)=\sigma(3,1)=8{\rm li}_4\frac12-2\z(2)\ln^22+7\z(3)\ln2+\frac1{3}\ln^42-{53\over 8}\z(4).$$ Hence (43) furnishes :$$\sigma(2,2)=-8{\rm
li}_4\frac12+2\z(2)\ln^22-\frac1{3}\ln^42-7\z(3)\ln 2+\frac{151}{16}\z(4).$$

 Setting $s=4$  in (36.1) we
get$$\sigma(3,2)+2\sigma(2,3)=-\sigma(4,1)+4\sum_{p\geq 1}{H_p\over
(2p+1)^4}+8\lambda(4)\ln2. $$ By the previous results, this gives$$\sigma(3,2)=-{31\over 2}\z(5)+\frac{35}4\z(2)\z(3)=-16\lambda(5)+\frac{40}3\lambda(2)\lambda(3).$$
\vfill\eject

\noindent {\bf (4) Application of a result of [6] to a formula for $\sigma(2a-1,2)$ }

 In [6], the authors gives a formula which leads to the calculation of $\sigma(2a-1,2).$ They define, for $q\geq 2$,$$E_{p,q}=\sum_{n\geq 1}{1\over
n^q}\sum_{1\leq k\leq 2n}{1\over k^p}\cdot \eqno(44)$$

With our notations, this turns to$$E_{p,q}=\sigma(q,p)+{1\over2^p}\z^*(q,p),$$ where $$\z^*(q,p)=\sum_{1\leq m\leq n}{1\over n^qm^p}=\z(p+q)+\z(q,p).$$ In particular, we get
$$\sigma(2a-1,2)=E_{2,2a-1}-\frac14\z^*(2a-1,2)$$ Theorem 2 in [6] give a formula for $E_{2,2a-1}$ which may be written as
:$$\eqalign{E_{2,2a-1}=\sum_{j=1}^{a-2}j2^{2j}\z(2j+1)\z(2a-2j)&+\Big((2a+1)2^{2a-3}-\frac12\Big)\z(2)\z(2a-1)\cr&-\Big(a2^{2a-1}+{2a^2-a-1\over
8}\Big)\z(2a+1).}\eqno(45)$$ But (see [1] or [2], th. (3.1)) $$\z^*(2a-1,2)=-\frac12(2a^2+a-1)\z(2a+1\  )+(2a-1)\z(2)\z(2a-1)+2\sum_{j=1}^{a-2}j\z(2j+1)\z(2a-2j),$$ so we
obtain, by some manipulations:$$\tableau{\eqalign{\sigma(2a-1,2)=-a2^{2a-1}\lambda(2a+1)+{2^{2a-1}(2a+1)\over
3}\lambda(2)\lambda(2a-1)\cr+\sum_{j=1}^{a-2}j2^{2j}\lambda(2j+1)\z(2a-2j)}\crh}.\eqno(46)=(g)$$ When $a=2$ this last formula gives again
$$\sigma(3,2)=-16\lambda(5)+\frac{40}3\lambda(2)\lambda(3).$$ For $a=3$, (46) gives
:$$\sigma(5,2)=-96\lambda(7)+\frac{224}3\lambda(2)\lambda(5)+4\lambda(3)\z(4)=-96\lambda(7)+\frac{224}3\lambda(2)\lambda(5)+{64\over15}\lambda(3)\lambda(4)$$
\vskip 2mm

\noindent {\bf (5) A  first weighted sum of $\sigma $'s .}  

Generally, the relation (36.1) allows us to obtain the sum $\ds\sum_{i=1}^{s-2}2^{i-1}\sigma(s-i,1+i)$ when $s$ is even and $\geq 4$ . Let be $s=2a, a\geq 2$ 
$$\sum_{i=1}^{2a-2}2^{i-1}\sigma(2a-i,1+i)=-J(2a)+2^{2a-2}h_{2a}+2^{2a-1}\lambda(2a)\ln2.\eqno(47)$$

In (47) $J(2a)$ is given by (23) and $h_{2a}$ by
(32.1):$$\eqalign{&\sum_{i=1}^{2a-2}2^{i-1}\sigma(2a-i,1+i)=-(2^{2a-1}\lambda(2a+1)-\sum_{j=1}^{a-1}2^{2j}\lambda(2j+1)\z(2a-2j))\cr
&+2^{2a-2}(-2\lambda(2a)\ln2+2a\lambda(2a+1)-2\sum_{j=1}^{a-1}\lambda(2j)\lambda(2a+1-2j))  +2^{2a-1}\lambda(2a) \ln 2.}\eqno (48)$$

In (47) the term in $\lambda(2a) \ln 2$
vanishes:$$\eqalign{\sum_{i=1}^{2a-2}2^{i-1}\sigma(2a-i,1+i)=&-(2^{2a-1}\lambda(2a+1)-\sum_{j=1}^{a-1}2^{2j}\lambda(2j+1)\z(2a-2j))\cr
&+2^{2a-2}(2a\lambda(2a+1)-2\sum_{j=1}^{a-1}\lambda(2j)\lambda(2a+1-2j))\cr =&
(a-1)2^{2a-1}\lambda(2a+1)+\sum_{j=1}^{a-1}2^{2a-2j}\z(2j)\lambda(2a+1-2j)\cr&-2^{2a-1}\sum_{j=1}^{a-1}\lambda(2j)\lambda(2a+1-2j) .}$$ In terms of series $\lambda$
this is :$$\sum_{i=1}^{2a-2}2^{i-1}\sigma(2a-i,1+i)=2^{2a-1}[
(a-1)\lambda(2a+1)+\sum_{j=1}^{a-1}{3-2^{2j}\over2^{2j}-1}\lambda(2j)\lambda(2a+1-2j)], $$
but it is lightly simpler to write$$\sum_{i=1}^{2a-2}2^{i-1}\sigma(2a-i,1+i)=2^{2a-1}[
(a-1)\lambda(2a+1)+\sum_{j=1}^{a-1}(3\cdot 2^{-2j}-1)\z(2j)\lambda(2a+1-2j)]. $$

In the left hand member we pick the term corresponding to $i=2a-2$. So, in view of (42) we 
obtain

$$\eqalign{\sum_{i=1}^{2a-3}2^{i-1}\sigma(2a-i,1+i)=&-2^{2a-2}[(2a^2-3a+2)\lambda(2a+1)\cr&-
2\sum_{j=1}^{a-1}\Big({3-2^{2j}\over2^{2j}-1}+2a-2j\Big)\lambda(2j)\lambda(2a+1-2j).]}\eqno(49)$$
 
We return to (37), in which we put now
$t=2r=2$:$$\sum_{i=1}^{s-2}2^{i-1}(i+1)\sigma(s-i,2+i)=-2^{s-1}\lambda(s)\lambda(2)+2^{s-2}sh_{s+1}+2^{s-1}s\lambda(s+1)\ln2.$$

If $s=4$,  this gives $$2\sigma(3,3)+6\sigma(2,4)=-8\lambda(4)\lambda(2)+16h_5+32\ln2\lambda(5).$$
Taking account of the value of $h_5$ given after $(39.2)$ and the relation $\ds \lambda(4)\lambda(2)=\frac54\lambda(6)$ we
obtain$$\sigma(3,3)+3\sigma(2,4)=15\lambda(6)-8\lambda^2(3).\eqno(50)$$One can verify that making $s=t=3$  or $s=4,t=2$ in (36) gives the same result.
\vskip 8mm

 \noindent {\bf (6) Other linear links between the sum $\sigma(k,\ell )$}\vskip 2mm

Now we give new linear relations  similar to those linking the $\z(k,\ell).$

For integers $k,\ell \geq 2$ we can write$$\lambda(k)\lambda(\ell)=\sum_{1\leq n\leq m}{1\over(2n-1)^k(2m-2n+1)^{\ell}}\cdot \eqno(51)$$

By using the decomposition in partial rational fractions given by (35), we get :

$$\eqalign{{1\over (2n-1)^k(2m-2n+1)^{\ell}}=&{1\over 2^{\ell}}\sum_{i=0}^{k-1}{\ell+i-1\choose i}{2^{-i}\over
m^{\ell+i}(2n-1)^{k-i}}\cr &+{1\over 2^k}\sum_{j=0}^{\ell-1}{k+j-1\choose j}{2^{-j}\over
m^{k+j}(2m-2n+1)^{\ell-j}}\cdot}\eqno (52)$$

So, by summation, (52) gives:$$\tableau{\eqalign{\lambda(k)\lambda(\ell)&={1\over 2^{\ell}}\sum_{i=0}^{k-1}2^{-i}{\ell+i-1\choose
i}\sigma(\ell+i,k-i)\cr&+{1\over 2^k}\sum_{j=0}^{\ell-1}2^{-j}{k+j-1\choose j}\sigma(k+j,\ell-j)}\crh}\,.\eqno (53)$$

In order to verify, test with $k=\ell=2$ : $4\lambda(2)^2=2(\sigma(2,2)+\sigma(3,1)),  $ which gives the
following  already seen relation:$$\sigma(2,2)+\sigma(3,1)=2\lambda(2)^2={\pi^4\over 32}\cdot$$
\vskip 3mm
\noindent{\bf VI. Sum formula for the $\sigma$'s of same weight}

\vskip 2mm We observe that $J(2)=\sigma(2,1)=2\lambda (3)$ , by (23), while  the relation (46) gives us $\sigma(2,2)+\sigma(3,1)=3\lambda(4)$. By
using the previous relations giving $\sigma(4,1), \sigma(3,2)$ and $\sigma(2,3)$, it is easy to see that the sum of the $\sigma$ series of weigth 5 is $4\lambda (5).$

So, this suggests to us that the following formula  $\ds\sum_{i=1}^{w-2}\sigma(w-i,i)=(w-1)\lambda(w)$ holds  for all weigths $w\geq 3$. Particular checking of this
relation for $w\leq 10$ have given to us the idea of the proof of the general result. This result is similar to those concerning the classical $MZV$: the so-call
{\sl theorem of the sum} asserts that the sum of all $MZV$ of same depth (number of variable) and of same weight $w$ is $\z(w)$. For a proof, see Zagier [7] or
Granville [3]. When the depth is 2, this gives $\ds\sum_{i=1}^{w-2}\z(w-i,i)=\z(w)$ (note that this last relation appears naturally in some proofs of Euler
formula (2.1)). Since $\ds \z^*(s,t)=\sum_{1\leq m\leq n}{1\over n^sm^t}=\z(s,t)+\z(s+t)$  we have: 
$\ds\sum_{i=1}^{w-2}\z^*(w-i,i)=(w-1)\z(w),$ which gives us a further motivation for the following theorem.

\proclaim Theorem. [of the $\sigma $-sum] The sum of all the $\sigma$'s of same weight is calculable: for $w\geq 3$, we have:
$$\tableau{\ds\sum_{i=1}^{w-2}\sigma(w-i,i)=(w-1)\lambda(w)\crh}.\eqno (54)=(h)$$

Since $\sigma(2,1)=2\lambda (3)$(see section V.3.3), it suffices to prove (54) when $w\geq 4$.
We start by transform the relation (53), for $w=k+\ell$. First we set $$x_i=\sigma(w-i,i)$$ for $1\leq i\leq w-2$. Put  $i'=k-i$ in the first summation of (53). So
it becomes :$${1\over 2^{\ell}}\sum_{i'=1}^{k}2^{i'-k}{w-i'-1\choose k-i'}\sigma(w-i',i')={1\over 2^w}\sum_{i=1}^k 2^{i}{w-i-1\choose \ell-1}x_i={1\over
2^w}\sum_{i=1}^{w-2} 2^{i}{w-i-1\choose \ell-1}x_i, $$ because the binomial numbers in the last sum vanish if $k+1\leq i\leq w-2. $ By procceding in same way for the
second sum in (53) we obtain finally:$$\lambda(k)\lambda(\ell)={1\over
2^w}\sum_{i=1}^{w-2} 2^{i}\bigg[{w-i-1\choose \ell-1}+{w-i-1\choose k-1}\bigg]x_i.\eqno (55)$$  

At this stage, the   proof depends on wether the weight is odd or not.\vskip 3mm

\noindent{\bf VI.1 The case $w$ odd}\vskip 2mm

We set $w=2a+1$. Return to the relation (42), which can be write:$$x_{2a-1}=2a(2a-1)\lambda(2a+1)-8A,\eqno (56.1)$$ by setting $\ds
A=\sum_{j=1}^{a-1}j\lambda(2a-2j)\lambda(2j+1).$

 But we have, by (55) applied to $k=2j+1$, $\ell=2a-2j, 1\leq j\leq a-1,$ $$2^{2a+1}A=\sum_{i=1}^{2a-1}2^{i}x_i\sum_{j=1}^{a-1}j\bigg[{2a-i\choose 2a-2j-1}+{2a-i\choose
2j}\bigg].\eqno (56.2)$$
 We denote by $A_i$ the inner sum of relation (56.2): we  prove that $2^{i}A_i$ does not depend of $i$ when $1\leq i\leq 2a-2$. 

In the first sum of $A_i$, we exchange $j$ into $a-j$. So we
get:$$A_i=a\sum_{j=1}^{a-1}{2a-i\choose2j-1}+\sum_{j=1}^{a-1}j{2a-i\choose2j}-\sum_{j=1}^{a-1}j{2a-i\choose2j-1}.\eqno (57)
$$Then we evaluate each  sum in (57). We use the binomial expansions$$(1+x)^{2a-i}+(1-x)^{2a-i}=2\sum_{0\leq j\leq a-\frac{i}2}{2a-i\choose
2j}x^{2j}=2\sum_{0\leq j\leq a-1}{2a-i\choose 2j}x^{2j},\eqno (58)$$  $$(1+x)^{2a-i}-(1-x)^{2a-i}=2\sum_{1\leq j\leq a-\frac{i-1}2}{2a-i\choose
2j-1}x^{2j-1}=2\sum_{1\leq j\leq a-1}{2a-i\choose 2j-1}x^{2j-1}\eqno (59)$$ and their
derivatives:$$(2a-i)((1+x)^{2a-i-1}-(1-x)^{2a-i-1})=4\sum_{0\leq j\leq a-1}j{2a-i\choose 2j}x^{2j-1},\eqno (60)$$
$$(2a-i)((1+x)^{2a-i-1}+(1-x)^{2a-i-1})=2\sum_{1\leq j\leq a-1}(2j-1){2a-i\choose 2j-1}x^{2j-2}.\eqno (61)$$

In the last formulas, we have introduce null binomial numbers in order to extend  the summation to all integer lying between 1 and $a-1$. 

Now, suppose that $1\leq i\leq 2a-2$. In (59), (60), (61) we set $x=1$, so this gives $$2\sum_{0\leq j\leq a-1}{2a-i\choose 2j-1}=2^{2a-i},\eqno(62)$$
$$4\sum_{1\leq j\leq a-1}j{2a-i\choose 2j}=(2a-i)2^{2a-i-1},\eqno (63)$$   $$2\sum_{1\leq j\leq a-1}(2j-1){2a-i\choose 2j-1}=(2a-i)2^{2a-i-1}.\eqno (64)$$ By relations (62) and
(64) we obtain:$$4\sum_{1\leq j\leq a-1}j{2a-i\choose 2j-1}=(2a-i)2^{2a-i-1}+2^{2a-i}.\eqno(65)$$ Then, we put this results in relation (57)( when
$1\leq i\leq 2a-2$):$$A_i=a2^{2a-i-1}-2^{2a-i-2}=2^{2a-2-i}(2a-1).$$So we get finally, when $1\leq i\leq 2a-2$:  $$2^{i}A_i=2^{2a-2}(2a-1).\eqno (66)$$ It remains the case when 
$i=2a-1$. Without difficulty , one obtains  :$$A_{2a-1}=a-1.\eqno(67)$$
Now we put  this calculations in (56.2). This gives $$2^{2a+1}A=2^{2a-2}(2a-1)+2^{2a-1}x_{2a-1}(a-1),$$ or$$8A=(2a-1)\sum_{i=1}^{2a-2}x_i+(2a-2)x_{2a-1}.$$ We put
this last result in (56.1) and the theorem is proved for odd weight.\vskip 3mm

\noindent{\bf VI.2  The case $w$ even}\vskip 2mm

The proof is similar, but easier. Let be $a\geq 2$. In the relation (55), we set $k=2j$ and $\ell =2a-2j$, for $1\leq j\leq a-1.$ So we have $w=2a$ and
$$2^{2a}\lambda(2j)\lambda(2a-2j)=\sum_{i=1}^{2a-2} 2^{i}\bigg[{2a-i-1\choose
2a-2j-1}+{2a-i-1\choose 2j-1}\bigg]x_i.\eqno (68)$$  Summing the last equalities for $1\leq j\leq a-1$ and using the formula $(E)$ of section V.2 we get
:$$2^{2a-1}(2a-1)\lambda (2a)=\sum_{i=1}^{2a-2}2^{i}x_i\sum_{j=1}^{a-1}\bigg[{2a-i-1\choose
2a-2j-1}+{2a-i-1\choose 2j-1}\bigg].\eqno (69)$$
But, by changing $j$ into $a-j$, we see that $\ds \sum_{j=1}^{a-1}{2a-i-1\choose
2a-2j-1}=\sum_{j=1}^{a-1}{2a-i-1\choose
2j-1}$. So it follows equality:$$2^{2a-1}(2a-1)\lambda (2a)=\sum_{i=1}^{2a-2}2^{i}x_i2\sum_{j=1}^{a-1}{2a-i-1\choose 2j-1}.\eqno (70)$$  As in the odd case, we set
$x=1$ in the binomial expansion :$$(1+x)^{2a-i-1}-(1-x)^{2a-i-1}=2\sum _{1\leq j\leq a-\frac{i}2}{2a-i-1\choose 2j-1}x^{2j-1}=2\sum _{1\leq j\leq a-1}{2a-i-1\choose
2j-1}x^{2j-1},$$ (if $\ds j>a-\frac{i}2,$ the corresponding binomial coefficient vanishes). So we get :$$2\sum _{1\leq j\leq a-1}{2a-i-1\choose
2j-1}=2^{2a-i-1}.\eqno (71)$$ By transporting this value in (70), we obtain$$\sum_{i=1}^{2a-2}x_i=(2a-1)\lambda(2a),$$ which is the wanted result for $w=2a$.\vskip 2mm

\noindent  {\bf VII.The explicit evaluation of other $\sigma$ series}\vskip 2mm

\noindent{\bf (1) The case of odd weights}

Suppose that $w$ is odd: the previons sections give explicit values   of the series $\sigma(2, w-2)$, $\sigma(w-2, 2)$ and 
$\sigma(w-1,1)$ in terms of
$\z$ or $\lambda $ series. When
$w=7$, the unknown series are
$ \sigma(4,3),\, \sigma(3,4)$ and the relation (55) fournishes easily the system  (where $x_i=\sigma(7-i,i)$)
$${\eqalign{6x_3+8x_4&=32\lambda(2)\lambda(5)-5x_1-5x_2-8x_5\cr 4x_3+2x_4&=16\lambda(3)\lambda(4)-5x_1-5x_2}}$$ which gives, taking account of the values of $x_1$, $x_2$
and $x_5$:
$${\eqalign{\sigma(4,3)&=120\lambda(7)-96\lambda(2)\lambda(5)\cr \sigma(3,4)&=-80\lambda(7)+8\lambda(3)\lambda(4)+{176\over 3}\lambda(2)\lambda(5).}}$$

When $w=2a+1\geq 9$, the number of unknown $\sigma$ series  is $w-5=2a-4$, but the relations (55) give only $a-1$ (which is $<2a-4$) equations. We think that these
equations are independant and that the other linear relations given by the formula (36) for various choices of $s$ and $t$  (with $s+t=2a+1$) depend of them.

However, we have find an alternative but tedious method to explicit   other $\sigma $ series  by integral representation of the finite sums $S_n^{(t)},$  similar
of these exposed in section {\bf II.2}. As example, we sketch out the case of $\sigma (2a-1,2)$, which will give an other proof of formula (46). This method uses 
trigonometric series which can be expressed by Bernoulli polynomials. It may be useful to compare  our method with the computations given in [6].

First, we observe that, for all integer $k\geq 1$, $$\int_0^{\pi}x\cos(2k-1)x\,{\rm d}x={-2\over(2k-1)^2}\virg\eqno (72)$$
and$$-2S_n^{(2)}=\int_0^{\pi}x\Big(\sum_{k=1}^n\cos(2k-1)x\Big)\,{\rm d}x=\frac12\int_0^{\pi}{x\sin 2nx\over \sin x}\,{\rm d}x. $$ It
follows:$$-4\sigma(2a-1,2)=\int_0^\pi{\varphi(x)\over \sin x}\,{\rm d}x,$$where$$\varphi(x)=x\sum_{n\geq 1}{\sin2nx\over n^{2a-1}}\cdot$$
As in in section {\bf II.2}, we set $\ds v_m=\int_0^\pi\varphi(x)\sin(2m-1)x\,{\rm d}x $ and prove that $$\sigma(2a-1,2)=-\frac12\sum_{m=1}^\infty v_m.$$
Then, we calculate $v_m$ by successive integrations by part, using at the last step the well known series (when $0<x<\pi)$:$$\sum_{n=1}^\infty{\cos 2nx\over
n^2}=\z(2)-x(\pi-x)\,\,\,{\rm and}\,\,\,\sum_{n=1}^\infty{\sin 2nx\over n}=\frac\pi2-x$$

This process can be extended as follow. From (72), we obtain $${-2\over(2k-1)^2}=\bigg[{x^2\over 2}\cos(2k-1)x\bigg]_0^\pi+\frac12(2k-1)\int_0^\pi x^2\sin(2k-1)x\,{\rm d}x,$$ which
gives $$\frac4{(2k-1)^3}=\frac{\pi^2}{2k-1}-\int_0^\pi x^2\sin(2k-1)x\,{\rm d}x$$ and$$4S_n^{(3)}=\pi^2S_n-\int_0^\pi x^2 \sum_{1\leq k\leq n}\sin(2k-1)x\,{\rm d}x=\pi^2S_n-\int_0^\pi x^2
{\sin^2nx\over \sin x}\,{\rm d}x.$$So:$$4\sigma(b,3)=\pi^2\sigma(b,1)-\int_0^\pi {x^2 \over \sin x}\sum_{n=1}^\infty{\sin^2 nx\over n^b}\,{\rm d}x,$$ which gives, after similar
calculations :
$${\eqalign{\sigma(2a-2,3)=&a(2a-1)2^{2a-3}\lambda(2a+1)\cr&-(a-1)(2a+3)2^{2a-4}\z(2)\lambda(2a-1)\cr &-\sum _{j=2}^{a-2} j(2j-1)2^{2j-2}\z(2a-2j)\lambda(2j+1)}}$$
Step by step, this process can give also $\sigma(2a-3,4)$ etc.

\noindent{\bf (2) The case of even weights}

 When the weight is even, it seems that no $\sigma $ series can be evaluate by classical functions, a part the case $w=4$  which 
needs the special value
$\ds {\rm li_4}\frac12\cdot$ If we attempt to explicit this $\sigma$'s by using the previous integral method, we fall on the trigonometric Clausen 
series $\displaystyle \sum_{n\geq 1}{\cos2nx\over n^{2a-1}}\virg$  and their derivatives (when $a>1)$. This series are not elementary, a part the case $a=1$.
\vskip 5mm

\noindent REFERENCES \vskip 3mm

\noindent[1] D. Borwein, J. M. Borwein, and R. Girgensohn, {\sl Explicit evaluation of Euler sums}, Proc. Edinburgh Math.Soc. 38  (1995), pp 277-294

\noindent[2] P. Flajolet, B. Salvy, {\sl Euler Sums and Contour Integral Representations}, Experiment Math, 7(1) 1998 pp.15-35 

\noindent[3] A. Granville,  {\sl A decomposition of Riemann zeta function}, Analytic Number Theory, London Math. Soc. Lecture Notes Ser., Vol. 247,
Cambridge University Press, Cambridge 1997, pp 95-101.

\noindent[4] P.F. Jordan , {\sl Infinite sums of psi  functions}, Bull.  Amer. Math. Soc. (79) 4, 1973 pp. 681-683

\noindent[5] R. Sitaramachandrarao , {\sl A formula of S. Ramanujan}, J. of Numbers Theory, (25) 1987  pp.1-19

\noindent[6] Yao Ling and Minking Eie, {\sl On recurrence relations for the extensions of
Euler sums}, Rocky Mountains J. Math., volume 38, number 1, 2008

\noindent[7] Zagier,{ \sl Multiple zeta values}, unpublished manuscript, Bonn 1995.

\end